\RequirePackage[l2tabu, orthodox]{nag}

\documentclass[12pt]{amsart}
\usepackage{fullpage,url,amssymb,enumerate,colonequals}
\usepackage{mathrsfs} 

\usepackage[OT2,T1]{fontenc}
\DeclareSymbolFont{cyrletters}{OT2}{wncyr}{m}{n}
\DeclareMathSymbol{\Sha}{\mathalpha}{cyrletters}{"58}

\usepackage{color}

\newcommand{\Aff}{\mathbb{A}}

\newcommand{\F}{\mathbb{F}}

\newcommand{\PP}{\mathbb{P}}
\newcommand{\Q}{\mathbb{Q}}

\newcommand{\Z}{\mathbb{Z}}

\newcommand{\Fbar}{{\overline{\F}}}

\newcommand{\calA}{\mathcal{A}}

\newcommand{\calX}{\mathcal{X}}

\newcommand{\calZ}{\mathcal{Z}}

\DeclareMathOperator{\End}{End}

\DeclareMathOperator{\Spec}{Spec}

\newcommand{\et}{{\operatorname{et}}}

\newcommand{\HH}{{\operatorname{H}}}

\newtheorem{theorem}{Theorem}[section]

\theoremstyle{definition}

\newtheorem*{hypothesisZ}{Hypothesis~Z}
\newtheorem*{hypothesisZprime}{Hypothesis~Z$'$}

\theoremstyle{remark}
\newtheorem{remark}[theorem]{Remark}

\makeatletter
\g@addto@macro\bfseries{\boldmath} 
\makeatother

\usepackage{microtype}

\usepackage[
	backref,
	pdfauthor={Bjorn Poonen}, 
]{hyperref}
\usepackage[alphabetic,backrefs,lite,nobysame]{amsrefs} 

\begin{document}

\title{Using zeta functions to factor polynomials over finite fields}
\author{Bjorn Poonen}
\thanks{This research was supported in part by National Science Foundation grant DMS-1601946 and Simons Foundation grants \#402472 (to Bjorn Poonen) and \#550033.}
\address{Department of Mathematics, Massachusetts Institute of Technology, Cambridge, MA 02139-4307, USA}
\email{poonen@math.mit.edu}
\urladdr{\url{http://math.mit.edu/~poonen/}}
\date{October 2, 2017}

\begin{abstract}
In 2005, Kayal suggested that Schoof's algorithm
for counting points on elliptic curves over finite fields
might yield an approach to factor polynomials over finite fields
in deterministic polynomial time.
We present an exposition of his idea and then explain details
of a generalization involving Pila's algorithm for abelian varieties.
\end{abstract}

\maketitle

\section{Introduction}\label{S:introduction}

Factoring univariate polynomials over finite fields
is a solved problem in practice.
Known algorithms are fast,
and are proved to run in polynomial time 
if granted access to a source of randomness.
But the theoretical question of whether there exists a 
\emph{deterministic} polynomial-time algorithm remains open.
See the survey articles \cite{Lenstra1982}, \cite{Lenstra1990}, 
and \cite{VonZurGathen-Panario2001};
the last of these contains a very extensive bibliography.

In 1985, Schoof gave a deterministic polynomial-time algorithm to compute 
the number of points on a given elliptic curve 
over a finite field \cite{Schoof1985}*{Section~3}.
At the Mathematisches Forschungsinstitut Oberwolfach in July 2005, 
Neeraj Kayal suggested a way to use Schoof's algorithm
to attempt to factor polynomials over finite fields
in deterministic polynomial time.
The author, who was present, responded 
that one could use higher genus curves
or higher-dimensional abelian varieties in place of elliptic curves, 
and that these heuristically had a greater chance of success.

It seems that the only written record of the ideas of Kayal and the author
before now
is the 2006 master's thesis of Amalaswintha Wolfsdorf \cite{Wolfsdorf2006}.
She describes Kayal's idea for elliptic curves in detail,
and writes a few sentences on the higher genus case
based on a November~25,~2005 email from the present author.
Our purpose is to present a brief exposition of Kayal's idea
and to explain details of the generalization,
which is Theorem~\ref{T:generalization of Kayal} in this article.

\section{Schoof's algorithm}\label{S:Schoof}

Understanding Kayal's idea requires some knowledge of 
Schoof's algorithm, which we now recall, in the special case 
of a prime field $\F_p$.

\begin{theorem}[\cite{Schoof1985}*{Section~3}]
There exists a deterministic polynomial-time algorithm that takes as input
a prime $p$ and a Weierstrass equation of an elliptic curve $E$ over $\F_p$,
and outputs $\#E(\F_p)$.
\end{theorem}

Polynomial-time means polynomial in the size of the input,
which is of order $\log p$.

\begin{proof}[Sketch of proof]
Hasse proved that $\#E(\F_p) = p-a+1$ for some $a \in \Z$
satisfying $|a| \le 2\sqrt{p}$.
If one can compute $a \bmod \ell$
for all primes $\ell \ne p$ up to some bound $L$,
then an effective Chinese remainder theorem lets one
compute $a \bmod \prod_{\ell < L} \ell$.
If $L$ is chosen as a sufficiently large constant multiple of $\log p$,
then $\prod_{\ell < L} \ell > 4 \sqrt{p}$,
so $a \bmod \prod_{\ell < L} \ell$ determines $a$.

The Frobenius endomorphism $F$ of $E$ satisfies $F^2-aF+p=0$ in $\End E$.
In particular, $F^2-aF+p$ acts as $0$ 
on the $\ell$-torsion subscheme $E[\ell]$,
and this condition uniquely determines $a \bmod \ell$.
It remains to explain how to compute with these objects.
First, $E[\ell]$ is $\Spec R$ for some $\F_p$-algebra $R$
defined by $O(1)$ explicit equations of degree polynomial in $\ell$,
and these equations can be computed from the group law on $E$;
from this, one can compute an explicit multiplication table
for $R$ with respect to an $\F_p$-basis.
(With a little more work, following Schoof, 
one can work even more explicitly by using division polynomials,
but this does not generalize as easily.)
The action of $F$ on $E[\ell]$ is given by the $p$th power map on $R$,
whose action on $\F_p$-algebra generators can be computed explicitly
by writing the exponent $p$ in binary 
and using repeated squaring and multiplication.
Similarly, the action of $p$ (or any smaller integer) on $E[\ell]$
can be computed by writing $p$ in binary and using repeated doubling 
and addition on $E$.
Combining these lets one compute the action of $F^2-aF+p$ on $E[\ell]$
in time bounded by $P(\ell,\log p)$ for some polynomial $P$.
For each $\ell$, try $a=0,1,\ldots,\ell-1$
until the value mod $\ell$ is found that makes $F^2-aF+p$ kill $E[\ell]$.
The total running time is at most
$\sum_{\ell < L} \sum_{a=0}^{\ell-1} P(\ell,\log p)$,
which is polynomial in $\log p$.
\end{proof}

\section{Kayal's factoring idea}\label{S:Kayal}

For simplicity, suppose that $p$ is a large prime
and suppose that we are given the product $f(t)=(t-r_1)(t-r_2) \in \F_p[t]$ 
for some unknown distinct $r_1,r_2 \in \F_p$.
Let $B = \F_p[t]/(f(t))$, so $B$ is secretly isomorphic to $\F_p \times \F_p$.
Elements of $B$ are represented by polynomials of degree $\le 1$ in $\F_p[t]$.
A Weierstrass equation over $B$ with discriminant in $B^\times$ 
defines an elliptic scheme $E$ over $B$,
which secretly specializes to two elliptic curves over $\F_p$, 
say $E_1$ and $E_2$.

What happens if we blithely run Schoof's algorithm on $E$, 
as if $B$ were $\F_p$?
If $\#E_1(\F_p)=\#E_2(\F_p)$, 
then there exists $a \in \Z$ such that $F^2-aF+p=0$ in $\End E$,
and the algorithm runs as usual,
and outputs the common value $\#E_1(\F_p) = E_2(\F_p)$,
but we learn nothing about the factorization of $f(t)$.
Now suppose instead that $\#E_1(\F_p) \ne \#E_2(\F_p)$.
Write $\#E_i(\F_p)=p-a_i+1$ for $i=1,2$, so $a_1 \ne a_2$.
Then for some $\ell$, we have $a_1 \not\equiv a_2 \pmod{\ell}$.
Thus, when we check integers $a$ to see if
$F^2-aF+p$ kills $E[\ell]$,
which amounts to certain elements of $B$ vanishing,
we instead find an integer $a_1$ for which 
these elements of $B$ vanish mod $t-r_1$ but do not all vanish mod $t-r_2$.
Hence we discover a nontrivial factor of $f(t)$.

Heuristically it is likely that $\#E_1(\F_p) \ne \#E_2(\F_p)$,
since there are about $4\sqrt{p}$ possible values for the order
of an elliptic curve over $\F_p$.
If we are unlucky enough to have chosen $E$ so that $\#E_1(\F_p) = E_2(\F_p)$,
we can try again with a different $E$,
or use the same linear polynomials as Weierstrass coefficients
while replacing $f(t)$ by $f(t+1)$.
We do not have a proof, however, that a deterministic sequence of
such trials will succeed after polynomially many attempts.

\begin{remark}
The same approach can be tried to factor
$f(t)=(t-r_1)\cdots(t-r_d)$ for distinct $r_1,\ldots,r_d \in \F_p$,
by induction on $d$.
If, using obvious notation,
$\#E_1(\F_p),\ldots,\#E_d(\F_p)$ are not all equal,
then Schoof's algorithm will find a nontrivial factor $g$ of $f$,
and then we can apply the inductive hypothesis to factor $g$ and $f/g$.
\end{remark}

\begin{remark}
\label{R:Berlekamp}
Berlekamp \cite{Berlekamp1970} showed that one can reduce the problem 
of factoring polynomials in $\F_q[t]$ for arbitrary prime powers $q$
to the problem of factoring polynomials in $\F_p[t]$ 
with distinct roots all in $\F_p$; 
see also \cite{Lenstra1982}*{Sections 3 and~4} for another exposition of this.
\end{remark}

\section{Pila's algorithm}\label{S:Pila}

To generalize Kayal's approach to abelian varieties,
we need Pila's generalization of Schoof's algorithm.

Let $A$ be a $g$-dimensional abelian variety over $\F_p$.
Let $F$ be the Frobenius endomorphism of $A$.
For each prime $\ell \ne p$, we may form the $\ell$-adic Tate module
$T_\ell A \colonequals \varprojlim_n A[\ell^n]$.
Let $P(t)$ be the characteristic polynomial of $F$ acting on $T_\ell A$,
so $\deg P = 2g$.
A priori the coefficients of $P$ are in $\Z_\ell$,
but in fact they lie in $\Z$ and are independent of the choice of $\ell$.
Knowledge of $P$ is equivalent to knowledge of the zeta function $Z_A$.

An abelian variety $A$ over $\F_p$ can be described explicitly
by giving a positive integer $N$ and a finite list of homogeneous polynomials 
in $\F_p[x_0,\ldots,x_N]$ whose common zero locus in $\PP^N$ is $A$,
together with the addition morphism $A \times A \to A$
(and also the inversion morphism if desired)
in terms of explicit polynomial mappings on affine patches.
Pila's algorithm accepts such data as input, and outputs $P(t) \in \Z[t]$.
Its running time is bounded by a polynomial in $\log p$
whose degree and coefficients depend only
on $N$ and the number and degrees of the polynomials
defining $A$ and the addition law \cite{Pila1990}*{Theorem~A}.

The general outline of Pila's algorithm is similar to that of Schoof's
algorithm: it computes $P(t) \bmod \ell$ for many small primes $\ell$
by studying the action of $F$ on $A[\ell]$, 
and then reconstructs $P(t)$ by using an effective Chinese remainder theorem.
For each $\ell$, it tries each monic degree~$2g$ polynomial 
in $(\Z/\ell \Z)[t]$ and tests whether it equals $P(t) \bmod \ell$.
Each test involves a deterministic sequence of 
arithmetic operations on elements of $\F_p$
controlled by queries:
each query asks whether some previously computed element is $0$,
and the result dictates which arithmetic operation is to be carried
out next.
This is all that we will need to know about Pila's algorithm.

\section{Generalization of Kayal's factoring idea}\label{S:generalization}

We will prove that we can replace elliptic curves 
by abelian varieties in Kayal's approach.
Also, instead of using only the order of the group of points, 
we can use the whole zeta function.
The advantage of using higher-dimensional abelian varieties 
is that there are many more possible zeta functions,
so success becomes very likely, at least heuristically: 
see Section~\ref{S:heuristic}.

Given a variety $V$ over a finite field, 
let $Z_V$ be its zeta function, viewed as a rational function in $\Q(T)$.

Let $U$ be a dense open subscheme of $\Aff^1_{\Z} \colonequals \Spec \Z[t]$.
Let $\calA \to U$ be an abelian scheme.
For each prime $p$ and $u \in U(\F_p)$,
let $\calA_u$ be the fiber above $u$.
Concretely, $\calA$ can be thought of as a family of abelian varieties
defined by equations with coefficients in $\Z[t]$; 
specializing $t$ to a suitably general value $u \in \F_p$ 
produces an abelian variety $\calA_u$ over $\F_p$;
here ``suitably general'' means outside a certain bad locus,
which may be taken to be of the form $\Delta(t)=0$
for some ``discriminant'' $\Delta(t) \in \Z[t]$ that is not identically zero
but vanishes at any $u$ for which the specialization $\calA_u$ is degenerate.

Theorem~\ref{T:generalization of Kayal} below will involve the following:

\begin{hypothesisZ}
There exist a dense open subscheme $U$ of $\Aff^1_\Z$
and an abelian scheme $\calA \to U$
such that for every sufficiently large prime $p$,
the $Z_{\calA_u}$ for the different $u \in U(\F_p)$ are distinct.
\end{hypothesisZ}

\begin{theorem}
\label{T:generalization of Kayal}
There is a deterministic algorithm that takes as input 
a finite field $\F_q$ and a nonzero polynomial $f \in \F_q[t]$,
and outputs the factors of $f$ in $\F_q[t]$,
such that if Hypothesis~Z holds,
then the running time is polynomial in $\log q$ and $\deg f$.
\end{theorem}

\begin{remark}
Our proof will show that 
an algorithm as in Theorem~\ref{T:generalization of Kayal} not only exists, 
but also can be written down explicitly, 
even if we do not know in advance 
the abelian scheme $\calA \to U$ in Hypothesis~Z.
\end{remark}

\begin{proof}
By Remark~\ref{R:Berlekamp},
we may assume that $q$ is a prime $p$ 
and that $f$ has distinct roots all in $\F_p$.
Of course we also assume that $\deg f \ge 2$.

First, we give an algorithm depending on explicit knowledge
of an abelian scheme $\calA \to U$ as in Hypothesis~Z.
More precisely, we may shrink $U$ to assume that $U = \Spec T$,
where $T=\Z[t][1/\Delta]$ for some nonzero polynomial $\Delta \in \Z[t]$,
and we may assume that we are given explicit polynomials 
describing $\calA$ over $T$ 
in the same way that we described abelian varieties over $\F_p$
in Section~\ref{S:Pila}.

There are at most $(\deg f)(\deg \Delta)$ values $c \in \F_p$
such that $f(t+c)$ and $\Delta(t)$ have a nontrivial gcd,
so by trying $c=0,1,\ldots$ in turn, we quickly find such a $c$
(of course, we may assume that $p > (\deg f)(\deg \Delta)$.
Replace $f(t)$ by $f(t+c)$ to assume that $\gcd(f,\Delta)=1$.
Let $B = \F_p[t]/(f(t))$.
Then $\Spec B$ is a closed subscheme of $U$,
and the base change $\calA_B$ is an abelian scheme over $B$.
It consists of a disjoint union of abelian varieties $\calA_u$ over $\F_p$,
one for each zero $u$ of $f$.

Apply Pila's algorithm to $\calA_B$,
but each time it queries an element of $B$ to test whether it is $0$,
instead compute a gcd with $f$ to test whether it is $0$, a unit,
or a nonzero zerodivisor.
By Hypothesis~Z, the zeta functions $Z_{\calA_u}$ 
for two different zeros $u$ of $f$ are distinct in $\Q(T)$,
so there exists a prime $\ell$ such that the characteristic polynomials
mod $\ell$ of $\calA_u$ for these two $u$ are distinct.
Therefore the computations in Pila's algorithm must eventually diverge
for these two values of $u$,
which can happen only if a nonzero zerodivisor in $B$ is encountered.
At that point, we have found a nontrivial factor $f_0$ of $f$.
Apply induction to the factors $f_0$ and $f/f_0$.
This completes the description of the algorithm 
when we are given $\calA \to U$ explicitly.
In particular, there exists an algorithm to factor polynomials,
even though we might not know which algorithm it is that does it.

We now describe a new program $\Omega$ that does not rely
on knowledge of $\calA \to U$.
Program $\Omega$ enumerates all computer programs
and runs them in parallel, 
devoting a fraction $2^{-n}$ of its computing power to the $n$th program;
at each step of each program, $\Omega$ tests whether what 
that program has printed so far 
is a list of linear polynomials over $\F_p$ whose product is $f$,
and if so, $\Omega$ terminates the whole computation with this answer.
If Hypothesis~Z is true,
and $n$ is the number of the program 
described in earlier paragraphs 
using an abelian scheme $\calA \to U$ as in Hypothesis~Z,
then $\Omega$ finds the factorization 
in time bounded by $2^n$ times a polynomial,
but $2^n$ is a constant, so this is still polynomial 
in the size of the input.
If Hypothesis~Z is false,
then $\Omega$ still terminates with the correct factorization 
because there exists $N$ such that program $N$ factors polynomials
by trial division, but the running time of $\Omega$ is not guaranteed
to be bounded by a polynomial in this case.
\end{proof}

\section{A heuristic for Hypothesis~Z}\label{S:heuristic}

For a $g$-dimensional abelian variety $A$ over $\F_p$,
the complex zeros of the characteristic polynomial $P(t)$ 
have absolute value $p^{1/2}$,
so the coefficient of $t^{2g-m}$ in $P(t)$
is $O_g(p^{m/2})$, with the implied constant depending on $g$ but not $p$.
Also, the functional equation of $Z_A$
shows that the coefficient of $t^m$ in $P(t)$ 
is determined by the coefficient of $t^{2g-m}$.
Thus $P(t)$ is determined by coefficients of $t^{2g-m}$
for $m=1,2,\ldots,g$, so there are at most 
$\prod_{m=1}^g O_g(p^{m/2}) = O_g(p^{g(g+1)/4})$ 
possibilities for $P(t)$.
Equivalently, if $\calZ_{g,p}$ is the set of zeta functions 
of all $g$-dimensional abelian varieties over $\F_p$,
then $\# \calZ_{g,p} = O_g(p^{g(g+1)/4})$ as $p \to \infty$.
In fact, DiPippo and Howe \cite{DiPippo-Howe1998} 
prove that for fixed $g$, 
we have $\# \calZ_{g,p} \sim p^{g(g+1)/4}$ as $p \to \infty$,
where in this section we use the notation $f(p) \sim h(p)$ 
to mean that $f(p)/h(p)$ tends to a positive constant 
depending only on $g$ as $p \to \infty$.

If we sample about $p$ zeta functions from $\calZ_{g,p}$ at random, 
then the expected number of equal pairs is 
$\sim \binom{p}{2} \frac{1}{p^{g(g+1)/4}} \sim p^{2-g(g+1)/4}$.
If we do this for all primes $p$ greater than or equal to some
large integer $p_0$,
then the expected total number of equal pairs for all $p$ is
$\sum_{p \ge p_0} p^{2-g(g+1)/4}$,
which tends to $0$ as $p_0 \to \infty$, 
provided that $2-g(g+1)/4 < -1$, which holds for $g \ge 4$.

Now let $U$ be a dense open subscheme of $\Aff^1_\Z$,
and let $\calA \to U$ be an abelian scheme of relative dimension~$g$.
The previous paragraph suggests that 
if we model the zeta functions of the fibers of $\calA \to U$
above $\F_p$-points of $U$ 
as being independent random elements of $\calZ_{g,p}$,
then for sufficiently large $p_0$
it is true for every $p \ge p_0$ that these zeta functions 
will be distinct;
in other words, $\calA \to U$ should satisfy the condition in Hypothesis~Z,
unless there is some extra structure to the family that the model
fails to reflect.

It even seems reasonable to guess that 
for a typical $1$-parameter family of genus~$4$ curves, 
the family of Jacobians will satisfy the condition in Hypothesis~Z.
See \cite{Sutherland-Voloch-preprint} for some specific candidate families.

\section{Weakening Hypothesis~Z}\label{S:weakening}

Something slightly weaker than Hypothesis~Z would suffice
to obtain a polynomial running time 
in Theorem~\ref{T:generalization of Kayal}:

\begin{hypothesisZprime}
There exist a dense open subscheme $U$ of $\Aff^1_\Z$ and
an abelian scheme $\calA \to U$ 
such that for every sufficiently large prime $p$,
there are at least $p-(\log p)^{O(1)}$ distinct zeta functions $Z_{\calA_u}$ 
as $u$ varies over $U(\F_p)$.
\end{hypothesisZprime}

Under Hypothesis~Z$'$, given $f \in \F_p[t]$ that factors completely,
one can attempt to factor the polynomials $f(t+c)$ 
for $(\log p)^{O(1)}$ different values $c \in \F_p$
by running Pila's algorithm 
as in the proof of Theorem~\ref{T:generalization of Kayal}.
As long as the $O(1)$ here is larger than the $O(1)$ in Hypothesis~Z$'$,
and as long as $p$ is sufficiently large,
there will be at least one such $c$
such that all the zeros of $f(t+c)$ mod $p$ lie in $U(\F_p)$
and the zeta functions of the fibers above these zeros
are pairwise distinct.
Thus the algorithm will succeed in factoring $f(t+c)$ for at least one $c$,
and evaluating the factors at $t-c$ recovers the factorization of $f(t)$.
The running time of the algorithm is still polynomial in $\log p$,
albeit possibly with a larger exponent.

\section{Using varieties other than abelian varieties}

Suppose that instead of an abelian scheme as in Hypothesis~Z, 
one had an arbitrary finite-type morphism $\calX \to U$
for a dense open subscheme $U$ of $\Aff^1_\Z$
such that for any sufficiently large prime $p$
and distinct $u_1,u_2 \in U(\F_p)$,
there exists a prime $\ell$ bounded by a polynomial in $\log p$
and a nonnegative integer $i$
such that the characteristic polynomials of Frobenius
acting on $\HH^i_{\et}(\calX_{u_1} \times \Fbar_p,\Z/\ell\Z)$
and $\HH^i_{\et}(\calX_{u_2} \times \Fbar_p,\Z/\ell\Z)$
are different.
Then again one could factor polynomials over finite fields
in deterministic polynomial time,
provided that one had an analogue of Pila's algorithm
that could compute these characteristic polynomials
using a deterministic sequence of arithmetic operations and queries 
whose number is bounded by a polynomial in $\ell$ 
whose degree and coefficients depend only on $\calX \to U$.

Madore and Orgogozo \cite{Madore-Orgogozo2015}*{Th\'eor\`eme~0.1} 
gave an algorithm for computing such characteristic polynomials,
but their bound on the running time is only primitive recursive,
not polynomial in $\ell$.

\section*{Acknowledgements} 

I thank Andrew Sutherland and Jos\'e Felipe Voloch 
for encouraging me to write this article,
and I thank Amalaswintha Wolfsdorf for sharing her master's thesis with me.

\end{document}